\documentclass[a4paper,12pt]{amsart}
\usepackage{amscd}
\theoremstyle{plain}
\newtheorem{theorem}{Theorem}[section]

\newtheorem{lemma}[theorem]{Lemma}
\newtheorem{cor}[theorem]{Corollary}
\theoremstyle{definition}
\newtheorem{remark}{Remark}
\newtheorem{assum}{Assumption}

\begin{document}

\def\pr{\mathbb{P}}
\def\com{\mathbb{C}}
\def\zet{\mathbb{Z}}
\def\real{\mathbb{R}}
\def\field{\mathbb{K}}
\def\razio{\mathbb{Q}}
\def\nat{\mathbb{N}}
\def\vl{X_\lambda}

\newcommand{\findim}{\hfill $\Box$ \par \smallskip}

\title[Lines on Calabi-Yau threefolds]{Calabi-Yau complete intersections with infinitely many
lines}
\author{Marcello Bernardara}

\begin{abstract}
We give two new examples of families of Calabi-Yau complete intersection
threefolds whose generic element contains infinitely many lines.
We get some results about the normal bundles of these lines and the Hilbert scheme
of lines on the threefolds.
\end{abstract}
\maketitle

\section{Calabi-Yau complete intersections and lines on them}
Throughout the paper, CY is used instead of Calabi-Yau.\par
The Clemens conjecture originally states that on the generic quintic
threefold the number of rational curves in a fixed homology class
is finite. More generally, the conjecture is expected to hold also
for CY complete intersection threefolds in ordinary
projective spaces (see \cite{johnknut}). In particular, all lines
on a CY threefold lie in the same homology class, hence the
conjecture states that the number of lines on the generic such threefold
is finite.\par Moreover, the expected number of
lines on a generic CY complete intersection threefold can
be computed with algebraic geometric techniques such as Schubert
calculus in the Grassmannians.\par We get the same result about CY
manifolds in mirror symmetry: there is a way to predict correctly the number $n_d$
of rational curves of a given degree $d$ lying on the generic
CY threefold.\par
Recall that a CY threefold is a complex compact K\"ahler threefold
$X$ with trivial canonical bundle:
$$K_X \simeq \mathcal{O}_X.$$
We will call a complete intersection of type $(d_1, \ldots, d_k)$ a
threefold which is a complete intersection of $k$ hypersurfaces in
${\pr}^{k+3}$ of degrees $d_1, \ldots, d_k$ respectively.\par The
adjunction formula for a complete intersection of type $(d_1,
\ldots, d_k)$
$$K_X \cong {\mathcal O}_X (\sum_{i=1}^k d_i - k - 4)$$
allows to conclude that the only projective CY
threefolds that are complete intersections are of type $(5)$,
(the quintic threefold in ${\pr}^4$), $(3,3)$ and $(4,2)$ in
${\pr}^5$, $(3,2,2)$ in ${\pr}^6$ and $(2,2,2,2)$ in
${\pr}^7$.\par Using Schubert calculus, we have the following
results about the number of lines on the generic threefold:\par
$$\begin{array}{rcl}
(5) &2875 &\mbox{lines} \\
(3,3) &1053 &\mbox{lines} \\
(4,2) &1280 &\mbox{lines} \\
(3,2,2) &720 &\mbox{lines} \\
(2,2,2,2) &512 &\mbox{lines.}
\end{array}$$
These results agree with mirror symmetry predictions (see
\cite{coxkatz, nordfjord} for mirror symmetry techniques,
\cite{libteit} for the case (3,3)).\par The
genericity assumption in Clemens conjecture is crucial: in fact
we know examples of CY threefolds with infinitely many
lines. The simplest is the Fermat quintic threefold in ${\pr}^4$,
defined by the equation
$$x_0^5 + x_1^5 + x_2^5 + x_3 ^5 + x_4^5 = 0.$$
The lines on the threefold are described in \cite{albano2}.\par
The first nontrivial example is due to van
Geemen. He found infinitely many lines on the generic
threefold of a family called the Dwork pencil. Its equation is
$$x_0^5 + x_1^5 + x_2^5 + x_3^5 +x_4^5 -5 \lambda x_0 x_1 x_2 x_3 x_4 = 0$$
hence we get a pencil of quintic hypersurfaces in $\pr^4$, whose
zero fiber is the Fermat.\par To see how to find lines on them,
see \cite{albano1}, for a deeper investigation, see \cite{mustata}.
The result is obtained by showing that on the
generic threefold of the family there are more than the expcted 2875 lines.
This can be done choosing a "good" automorphism of the threefold and
finding lines fixed by it. In this case "good" means that its order does
not divide the expected number of lines, so it has fixed lines.
Under the action of the automorphisms of the threefold, the orbit of one of those
contains at least 5000 lines, clearly more than 2875.

\section{A $(3,3)$ complete intersection pencil}\label{33}

\begin{assum} We consider exclusively projective spaces over the
complex field $\com$.\end{assum}\par

The first example is a pencil of
CY threefolds of type $(3,3)$. On this particular pencil,
we are able to construct more lines than expected (1458 instead of
1053).\par The equations of the generic threefold
$X_{\lambda}$ (smooth for $\lambda$ generic) of the pencil are:
\begin{equation}\label{defpen}
\vl :=
\begin{cases}
x_0^3 + x_1^3 + x_2^3 - 3 \lambda x_3 x_4 x_5 = 0\\
x_3^3 + x_4^3 + x_5^3 - 3 \lambda x_0 x_1 x_2 = 0.
\end{cases}
\end{equation}
This pencil is invariant under a group of automorphisms
of ${\pr}^5$ of order 81 (see \cite{libteit}).\par Let $\phi$ be
the involution of $\pr^5$ given by the change of coordinates
$(12)(45) \in S_6$, which preserves $\vl$. We consider its invariant
subspaces $V_{\pm}$:
\begin{equation}\label{autospazi}
\begin{array}{rl}
V_+ &= \lbrace (a : a : b : c : c : d) \rbrace \\
V_- &= \lbrace (q : -q : 0 : p : -p : 0) \rbrace.
\end{array}
\end{equation}
Consider lines either contained in one of these subspaces or
intersecting both; such lines are $\phi$-invariant. In this
case there is no line lying on $\vl$ entirely contained in $V_{\pm}$,
but we have the following result.
\begin{lemma}\label{costrlines} On the generic threefold
$X_{\lambda}$ there are 36  lines connecting the invariant
subspaces (\ref{autospazi}), hence each one is fixed by $\phi$.
\end{lemma}
\begin{proof} It can be easily seen that there are no points in $V_{\pm}$ lying on
$X_{\lambda}$ if $d=0$ or $q=0$, so, without loss of generality,
consider:
$$V_+ = \lbrace (a:a:b:c:c:1) \vert a,b,c \in {\com} \rbrace $$
$$V_- = \lbrace (1:-1:0:p:-p:0) \vert p \in {\com} \rbrace.$$
Lines joining such points have parametric equations:
\begin{equation}\label{param}
(at + s : at - s : bt :  ct +ps : ct - ps : t)
\end{equation}
where $(s:t)$ ranges over ${\pr}^1$.\par Substituting the equation
(\ref{param}) in the equations (\ref{defpen}) of
$X_{\lambda}$, we obtain two cubic homogeneous polynomials in $s$, $t$.
The line belongs to the threefold if and only
if these polynomials vanish identically. It appears in the following cases:
$$a^3 = \dfrac{(2c^3 + 1) \lambda}{12c}$$
$$b = \dfrac{4ac}{\lambda^2}$$
$$p^2 = -\dfrac{2a}{\lambda}$$
and $c$ satisfies:
$$64 c^6 - (16 \lambda^6 - 32)c^3 +\lambda^6 = 0.$$
In particular, we have 6 values for $c$ for generic $\lambda$,
then we have 18 values for $a$ and 36 for $p$.
\end{proof}
\begin{theorem}\label{rettesu33}
On the generic threefold in the pencil $X_{\lambda}$ there are infinitely many lines.
\end{theorem}
\begin{proof} We know from the preceding Lemma that we have 36 lines on
$\vl$. Pick one of them and call it $l$.\par Consider the action
of the group $({\com}^{\ast})^6$ on ${\pr}^5$, where an element
$(a_0,\ldots,a_5)$ $\in$ $({\com}^{\ast})^6$ acts componentwise by:
$$(a_0 ,a_1 ,a_2 ,a_3 ,a_4 ,a_5 ) \cdot (x_0 :x_1 :x_2 :x_3 :x_4 :x_5 ) =$$
$$(a_0 x_0 : a_1 x_1 : a_2 x_2 : a_3 x_3 : a_4 x_4 : a_5 x_5 )$$
on $(x_0 : \ldots : x_5 )$ $\in$ ${\pr}^5$.\par Let $\alpha_i$ in
$({\com}^{\ast})^6$ be the elements:
\begin{equation}\begin{array}{cc}
\alpha_1 &= (1, \omega, \omega^{-1}, 1,1,1)\\
\alpha_2 &= (1,1,\omega,\zeta,\zeta,\zeta^{-2})\\
\alpha_3 &= (1,1,1,1,\omega,\omega^{-1})
\end{array}
\end{equation}
where $\zeta$ is a primitive ninth root of unity and $\omega =
\zeta^3$.\par We note that the group
$$G := < \alpha_1, \alpha_2, \alpha_3 > \subset \mathrm{Aut} \vl,$$
has order 81.\par Two other subgroups of
$\mathrm{Aut} \vl$ are given by the actions of $S_3$
on the first three coordinates and on the last three.
We denote the product of these two groups by $H$.\par The orbit of
the line $l$ under the action of the group $G$ has order 81,
because no element of this group fixes $l$.\par The orbit of $l$
under the action of $H$ has 18 elements, because $\phi = (12)(45)$
fixes $l$.\par
Consider now the group $G \times H$ and check that if
$g h (l) = l$, where $g \in G$ and $h \in H$, then $g (l) = l$ and
$h (l) = l$ (consider the points of $h(l) \cap V_-$ and
then the action of $G$ on these points).\par Hence, the order of
the orbit of $l$ under the action of the group $G \times H$, is
$81 \cdot 18 = 1458$. This number is larger than expected.
\end{proof}
\begin{remark}\label{chernclass} Recall the
way of counting lines proposed by S. Katz in \cite{katz2}. It is
based on finding a compact moduli space $\mathcal{M}$ of the curves on the
manifold, then constructing a rank $r = \mathrm{dim} \mathcal{M}$ vector bundle
with some good properties and then computing its $r$-th Chern class. 
We note that, because of their construction with
automorphisms, our lines have the same behavior. If they
were isolated, each would count as one; this would make the
calculation fail. Then we deduce that each of these lines
belongs to a continuous family. Notice that this tells nothing
about the number and the geometric properties of these families,
except that this excludes the case that these lines have
normal bundle of the form ${\mathcal O}_{\pr^1} (-1)
\oplus {\mathcal O}_{\pr^1} (-1)$.\end{remark}\par

\section{Normal bundle and Hilbert scheme of lines}

Recall the definition of normal bundle of a line $L$ in a
manifold $X$, as the cokernel in the exact sequence:
\begin{equation}\label{definiznorm}0 \longrightarrow T_L \longrightarrow T_{X \vert
L} \longrightarrow N \longrightarrow 0.
\end{equation} Now we are looking
for the normal bundle of a line $L$ on a threefold $X$, which is a
bundle over ${\pr}^1$, hence we can split it as
\begin{equation}\label{aeb}
N \cong {\mathcal O}_{\pr^1} (a) \oplus {\mathcal O}_{\pr^1} (b).
\end{equation} From the CY condition we deduce
(see e. g. \cite{katz2})
$$a + b = -2.$$
Let $X$ be a projective variety and $Z \subset X$ a subvariety. It is
well known (see \cite{kollar}), that for the Zariski tangent space to
the Hilbert scheme in $\lbrack Z \rbrack$ the following isomorphism holds:
$$T_{\lbrack Z \rbrack} \mathrm{Hilb}(X) \cong \mathrm{Hom}_X (I(Z),{\mathcal O}_Z)
= \mathrm{Hom}_Z (I(Z)/I(Z)^2,{\mathcal O}_Z).$$ The right side is the 
zeroth cohomology group of the normal bundle of $Z$ in $X$ (see \cite{hartshorne}), thus:
\begin{equation}\label{isomorfismonormale}
T_{\lbrack Z \rbrack} \mathrm{Hilb}(X) \cong H^0 (N_{Z \vert X}).
\end{equation}
In our case, we are looking for the normal bundle of lines lying
in a continuous family, hence the Zariski tangent space to the Hilbert
scheme in the point corresponding to these lines should be
positive dimensional. This gives
$$N_{l \vert X_\lambda} \not\cong {\mathcal O}_{\pr^1}(-1) \oplus {\mathcal O}_{\pr^1}(-1)$$
because in this case we would have $h^0 (N) = 0$.\par

\subsection{How to calculate the normal bundle} In this section, we show how
to calculate the normal bundle of a line on a CY complete intersection threefold
and after we will apply the calculation to the lines previously constructed.\par
Our aim is to calculate $a$ and $b$ in (\ref{aeb}), trying to
generalize slightly the calculations in \cite{katz1} to the complete
intersection case.\par
First, let $X$ be a hypersurface in $\pr^n$ and $L \subset X$ a
line on it. Change the coordinates of ${\pr}^n$ such that the line $L$ has
parametrization $(s:t:0:\cdots:0)$; in this case, the ideal $I_L$ of $L$ is $I_L =
(x_2, \cdots, x_n)$. Let us call $F_d$ the polynomial defining
$X$ and $d$ its degree. $L \subset X$ and $L$ is
the intersection of the hyperplanes $x_2 = \cdots = x_n = 0$, so we
can write:
$$F_d = x_2 F_2 + \cdots x_n F_n.$$
Modulo elements of $I_L^2$, we get
$$F_d = x_2 f_2(x_0,x_1) + \cdots + x_n f_n (x_0,x_1)$$
where each $f_i$ is homogeneous of degree $d-1$; it can be seen as
$F_{i \vert L}$.\par Using these exact sequences
$$\begin{array}{cc}
a) &0 \longrightarrow N \longrightarrow {\mathcal
O}_{{\pr}^1} (1)^{n-1} \longrightarrow {\mathcal O}_{{\pr}^1}(d)
\longrightarrow 0\\
b) &0 \longrightarrow T_X \longrightarrow T_{{\pr}^n \vert X}
\longrightarrow {\mathcal O}_X (d) \longrightarrow 0 \end{array}
$$
and (\ref{definiznorm}), it is possible (\cite{katz1}) to get the
normal bundle as the kernel of the map
\begin{equation}\label{soll}\begin{array}{rl} {\mathcal
O}_{{\pr}^1} (1)^{n-1} &\longrightarrow {\mathcal
O}_{{\pr}^1}(d)\\
(s_2,\cdots,s_n) &\longmapsto \sum_{i=2}^n f_i s_i.
\end{array}
\end{equation}
Now let $X$ be a complete intersection of two hypersurfaces of
degree $d$ and $e$, given respectively by $F = 0$ and $G = 0$.
Let $L \subset X$ be parametrized as before, so we get the
homogeneous polynomials $f_i$ and $g_i$ of degrees $d-1$ and $e-1$
respectively in $(x_0, x_1)$.\par From a direct calculation, we
get
$$N_{L \vert X} \simeq ker ({\mathcal O}_{{\pr}^1}(1)^{n-1}
\stackrel{M}{\longrightarrow} {\mathcal O}_{{\pr}^1}(d) \oplus
{\mathcal O}_{{\pr}^1}(e))$$ where the map is given by the $2
\times (n-1)$ matrix $M$ with rows given by the $f_i$ and $g_i$.
If we let $A = {\com}\lbrack x_0, x_1 \rbrack$, we can rewrite
this map as a map $A^{n-1} \to A^2 $. Hence we are looking at the
module
$$B = ker(A^{n-1}(1) \stackrel{M}{\longrightarrow} A(d) \oplus A(e))$$
and we know (see for example \cite{hilb}) that $B$ has a basis of
vectors of homogeneous polynomials $T_i$ of the same degree (within
the vector) $t_i$. In the case of a line in a threefold we have $i
= 1,2$, hence:
$$N = {\mathcal O}_{\pr^1}(1-t_1) \oplus {\mathcal O}_{\pr^1}(1 - t_2).$$
In conclusion we get the following result.
\begin{theorem}\label{solnorm}
Let $T$ a vector of homogeneous polynomials of minimal degree $t$ in
$(x_0,x_1)$ such that
$$M \cdot T = 0$$
where $M$ is the matrix with rows given by the $f_i$ and the
$g_i$. Then the normal bundle $N_{L \vert X}$ splits in the following way:
$$\begin{array}{clc}
N_{L \vert X} &= {\mathcal O}_{\pr^1} (1) \oplus {\mathcal
O}_{\pr^1} (-3)
&\mbox{if $t = 0$}\\
N_{L \vert X} &= {\mathcal O}_{\pr^1} \oplus {\mathcal O}_{\pr^1}
(-2) &\mbox{if
$t = 1$}\\
N_{L \vert X} &= {\mathcal O}_{\pr^1} (-1) \oplus {\mathcal
O}_{\pr^1} (-1) &\mbox{otherwise.}\end{array}
$$
\end{theorem}
\begin{proof} This follows easily from the above considerations,
remembering that we should have, for the CY condition,
$t_1 - 1 + t_2 - 1 = -2$.
\end{proof}
The argument is essentially the same for a generic CY complete intersection
threefold in projective space.

\subsection{The normal bundle of constructed lines}
We calculate the normal bundle of the lines constructed in the previous
section on the generic threefold $X_\lambda$.\par
\begin{lemma}\label{calcolofibnorm} Let $\lambda$ be generic and $l \subset \vl$
be the line parametrized by
$$(at + s: at - s : bt : ct + ps : ct - ps : t)$$
as in Lemma \ref{costrlines}.\par Then its normal bundle on $\vl$ splits as:
$$N_{l \vert X_\lambda} \cong {\mathcal O}_{\pr^1} \oplus {\mathcal O}_{\pr^1}(-2).$$
\end{lemma}
\begin{proof} Recall that
$$N_{l \vert X_\lambda} \not\cong {\mathcal O}_{\pr^1} (-1) \oplus {\mathcal O}_{\pr^1}(-1).$$
Define new coordinates:
$$
\begin{cases}
x_0 &= y_0 + ay_5 \\
x_1 &= -y_0 + y_1 + ay_5 \\
x_2 &= y_2 + by_5 \\
x_3 &= py_0 + y_3 + cy_5 \\
x_4 &= -py_0 + y_4 + cy_5 \\
x_5 &= y_5.
\end{cases}
$$
$l$ has now parametrization $(s:0:0:0:0:t)$.\par We can obtain the matrix
$M$, with coefficients homogeneous quadratic polynomials in $y_0$ and $y_5$:
$$
\left(
\begin{matrix}
(y_0 - ay_5)^2 & b^2 y_5^2 & \lambda p(y_0 y_5) - c \lambda y_5^2 & -\lambda p(y_0 y_5) - c \lambda y_5^2 \\
-b \lambda (y_0 y_5 + y_5^2) & \lambda (y_0^2 - a^2 y_5^2) & (py_0 + cy_5)^2 & (py_0 - cy_5)^2
\end{matrix}
\right).
$$
Now we verify that there are no nonzero vectors $B \in {\com}^4$ such that $M
\cdot B = 0$. This leads to:
$$N_{l \vert X_\lambda} \not\cong {\mathcal O}_{\pr^1}(1) \oplus {\mathcal O}_{\pr^1}(-3).$$
\end{proof}
\begin{cor}
The same result holds for each line previously constructed on the generic
threefold $X_\lambda$.\end{cor}
\begin{proof} Let $l$ be as in Lemma \ref{calcolofibnorm}. Each line constructed 
in Theorem \ref{costrlines} can be obtained by $l$ using an automorphism of $\vl$.
\end{proof}\par
This lead us to conclude the dimension of the Hilbert scheme is
positive, in particular
$${\mathrm dim} T_{\lbrack l \rbrack} {\mathcal H}_\lambda = h^0 (N_{l \vert X_\lambda}) = 1$$
for the lines $l$ we constructed.\par

\section{A (2,2,2,2) two-parameter Family}

We now give a new example, a two-parameter family
of $(2,2,2,2)$ threefolds in ${\pr}^7$. Consider
the family (smooth for generic ($\lambda, \mu$)) obtained by the
complete intersection of the four quadrics
\begin{equation}\label{quad1}
X_{\lambda, \mu}:=
\begin{cases}
x_0^2 + x_1^2 + x_2^2 + x_3^2 + x_4^2 + x_5^2
\phantom{+ x_6^2 +x_7^2}- 2 \mu x_6 x_7 = 0\\
x_0^2 + x_1^2 + x_2^2 + x_3^2 \phantom{+ x_6^2 +x_7^2}+ x_6^2 +
x_7^2 - 2 \lambda x_4 x_5 = 0\\
x_0^2 + x_1^2 \phantom{+ x_6^2 +x_7^2}+ x_4^2 + x_5^2 + x_6^2 +
x_7^2 - 2 \lambda x_2 x_3 = 0\\
\phantom{+ x_6^2 +x_7^2}x_2^2 + x_3^2 + x_4^2 + x_5^2 + x_6^2 +
x_7^2 - 2 \lambda x_0 x_1 = 0.
\end{cases}
\end{equation}
As in the previous case, on the generic threefold more than the
512 expected lines are shown.\par The technique is the same: in this
case we take $\phi$ to be the order 3 automorphism of $X_{\lambda,
\mu}$ given by the permutation of coordinates $(135)(246)$ in $\pr^7$.
On $\pr^7$ we consider its invariant subspaces:
\begin{equation}
V_+ = \lbrace( a : b : a : b : a : b : c : d ) \rbrace \notag
\end{equation}
\begin{equation}
V_{\omega} = \lbrace( p : q : \omega p : \omega q : \omega^2 p :
\omega^2 q : 0 : 0 ) \rbrace \notag
\end{equation}
where $(a:b:c:d) \in {\pr}^3$, $(p:q) \in {\pr}^1$ and $\omega \in {\com}$
is primitive third root of the unity.
\begin{lemma}\label{rette2222}
For generic $(\lambda, \mu)$, there are 8 lines on $X_{\lambda, \mu}$
intersecting both $V_+$ and $V_\omega$.
\end{lemma}
\begin{proof} First we verify that the points with $b = 0$ and the ones with
$p = 0$ don't lie on the threefold, hence we can consider:
\begin{equation}\label{autos12}
V_+ = \lbrace( a : 1 : a : 1 : a : 1 : c : d ) \vert a,c,d \in
{\com} \rbrace
\end{equation}
\begin{equation}\label{autos22}
V_{\omega} = \lbrace( 1 : q : \omega : \omega q : \omega^2 :
\omega^2 q : 0 : 0 ) \vert q\in {\com} \rbrace.
\end{equation}
Lines joining these points have parametric equation
\begin{equation}\label{retta2}
(s + at : qs + t : \omega s + at : \omega qs + t : \omega^2 s + at
: \omega^2 qs + t : ct : dt)
\end{equation}
where $(s : t) \in {\pr}^1$.\par Substituting these values into the
equations of the generic threefold, the line lies on $X_{\lambda,
\mu}$ if and only if
\begin{align}
a &= - \dfrac{\lambda + q}{\lambda q +1} \notag \\
d &= \dfrac{3 a^2 + 3}{2 \mu c} \notag
\end{align}
where $q$ is a root of
$$q^2 + 2 \lambda q + 1 = 0$$
and $c$ is a root of
$$4 \mu^2 c^4 + ( 2 a^2 - 2 \lambda a + 2 ) c^2 + ( 3 a^2 + 3 )^2 = 0.$$
We get the proof, remarking that if we have 2 values for $a$, then it is clear
that we have 8 values for $c$ for the generic couple
$(\lambda, \mu)$ and that this does not depend on $a$ and on
$q$. \end{proof}
\begin{theorem}\label{teoquad}
On the generic threefold of the pencil $X_{\lambda, \mu}$ there
are infinitely many lines.
\end{theorem}
\begin{proof} A subgroup of Aut$X_{\lambda, \mu}$ is given by the action
of $S_3$ on the first three pairs of coordinates. The orbit of one of
the constructed lines under this automorphism group consists of 2 lines, because
$\phi$ belongs to this group.\par Another
subgroup of automorphisms is generated by the permutations of
coordinates $(12)$, $(34)$, $(56)$ and $(78)$: the constructed
lines are not fixed by any of these automorphisms, hence
the orbit of each line under the action of this subgroup has
16 elements.\par Let $G$ be the product of these two groups (in
particular, $G \leq S_8$).\par Let $H$ be the subgroup of
$({\com}^{\ast})^8$, acting on $\pr^7$ by the
coordinatewise product, with generators
$$\begin{array}{ccccccccc}
\alpha_1 = ( &-1 , &-1 , &1 , &1 , &1 , &1 , &1 , &1) \\
\alpha_2 = ( &1 , &1 , &-1 , &-1 , &1 , &1 , &1 , &1) \\
\alpha_3 = ( &1 , &1 , &1 , &1 , &-1 , &-1 , &1 , &1).
\end{array}$$
The orbit of each line under its action consists of 8
lines.\par Consider now the group $G \times H$: we have to
check that if $g h (l) = l$, where $g \in G$ and $h \in H$, then
$g (l) = l$ and $h (l) = l$ (consider the points in the set $h(l)
\cap V_{\omega}$ and then the action of $G$ on these points). We
get finally that the orbit of each line under the action of this
group consists of 256 elements.\par The key remark now is that the
orbits of the lines are disjoint, and this is made making a table
comparing the values obtained for the points in $V_\omega$ and $V_+$
starting from different values of $c$.\par We finally
get on $X_{\lambda,\mu}$ at least 2048 lines, that is more than
expected.\end{proof}
\begin{remark} The
same argument used in Remark \ref{chernclass}
is valid in this case, so all the constructed lines belong to a
continuous family.\end{remark}\par
\subsection{The normal bundle}
\begin{lemma}
The line $l \subset X_{\lambda, \mu}$ parametrized by
$$(s + at : qs + t : \omega s + at : \omega qs + t : \omega^2 s +
at : \omega^2 qs + t : ct : dt)$$ as in Lemma \ref{rette2222}, has
normal bundle
$$N_{l \vert X_{\lambda, \mu}} \cong {\mathcal O}_{\pr^1} \oplus {\mathcal O}_{\pr^1} (-2).$$
\end{lemma}
\begin{proof} As before:
$$N_{l \vert X_{\lambda, \mu}}
\not\cong {\mathcal O}_{\pr^1}(-1) \oplus {\mathcal
O}_{\pr^1}(-1).$$
Now with the change of coordinates:
$$
\begin{cases}
x_0 &= y_0 + a y_7 \\
x_1 &= qy_0 + y_1 + y_7 \\
x_2 &= \omega y_0 + y_2 + a y_7 \\
x_3 &= \omega qy_0 + y_3 + y_7 \\
x_4 &= \omega^2 y_0 + y_4 + a y_7 \\
x_5 &= \omega^2 qy_0 + y_5 + y_7 \\
x_6 &= y_6 + c y_7 \\
x_7 &= d y_7
\end{cases}
$$
the line $l$ gets parametrization $(s:0:0:0:0:0:0:t)$. We calculate the matrix $M$,
with coefficients linear homogeneous polynomials in $y_0$ and $y_7$: \tiny
$$
\left(
\begin{matrix}
q y_0 + y_7 & \omega y_0 + a y_7 & \omega q y_0 + y_7 & \omega^2
y_0 + a y_7 & \omega^2 q y_0 + y_7 & -\mu d y_7 \\
q y_0 + y_7 & \omega y_0 + a y_7 & \omega q y_0 + y_7 & -\lambda
(\omega^2 q y_0 + y_7) &  -\lambda (\omega^2 y_0 + a y_7) & c y_7 \\
q y_0 + y_7 & -\lambda ( \omega q y_0 + y_7) & -\lambda ( \omega
y_0 + a y_7) & \omega^2 y_0 + a y_7 & \omega^2 q y_0 + y_7 & c y_7\\
-\lambda (y_0 + a y_7) &  \omega y_0 + a y_7 & \omega q y_0 + y_7
& \omega^2 y_0 + a y_7 & \omega^2 q y_0 + y_7 & c y_7 \\
\end{matrix}
\right).
$$
\normalsize We now verify that for generic $(\lambda, \mu)$ there are
no nonzero vectors $B$ in ${\com}^6$
such that $M \cdot B = 0$ and then $N_{l \vert X_{\lambda, \mu}}
\not\cong {\mathcal O}_{\pr^1} (1) \oplus {\mathcal O}_{\pr^1}
(-3)$.
\end{proof}
We deduce that all constructed lines have such a normal bundle and
that the dimension of the Hilbert scheme is positive, in particular:
$${\mathrm dim} T_{\lbrack l \rbrack} {\mathcal H}_\lambda = h^0 (N_{l \vert X_(\lambda,\mu)}) = 1$$
for the lines $l$ we constructed.\par

Marcello Bernardara: Laboratoire J. A. Dieudonn\'e, Universit\'e de Nice Sophia Antipolis, 06108 Nice Cedex 2.\par
E-mail: bernamar@math.unice.fr
\end{document}